\documentclass[12pt,leqno]{amsart}
\textheight9in  \def\DATE{\today}
\newtheorem{theorem}{Theorem}
\newtheorem{definition}[theorem]{Definition}

\newtheorem{proposition}[theorem]{Proposition}

\newcommand\pa{3\mbox{\rm{-}}pa\mathcal{A}ss}
\newcommand\ta{3\mbox{\rm{-}}tot\mathcal{A}ss}

\newcommand\pd{3\mbox{\rm{-}}pa^1\mathcal{A}ss}
\newcommand\npa{n\mbox{\rm{-}}pa\mathcal{A}ss}
\newcommand\nta{n\mbox{\rm{-}}tot\mathcal{A}ss}
\newcommand\ntad{n\mbox{\rm{-}}tot^1\mathcal{A}ss}
\newcommand\tota{3\mbox{\rm{-}}tot_{\tau_{13}}\mathcal{A}ss}
\newcommand\pad{3\mbox{\rm{-}}pa_{\tau_{13}}^1\mathcal{A}ss}
\newcommand\pp{\noindent{\it Proof. }}

\title{The $n$-ary algebra of tensors and  of cubic and hypercubic matrices}

\author{Nicolas GOZE - Elisabeth REMM }
\thanks{LMIA. Universit\'e de Haute Alsace, 4 rue des Fr\`{e}res Lumi\`{e}re, F. 68093 Mulhouse Cedex\\
nicolas.goze@uha.fr, elisabeth.remm@uha.fr.}

\begin{document}

\maketitle

\begin{abstract}
We define a ternary product and more generally a $(2k+1)$-ary product on the vector space $T^p_q(E)$  
of tensors of type $(p, q)$ that is  contravariant of  order $p$, covariant of order $q$ and total order $(p+q)$. 
This product is totally associative up 
to a permutation $s_k$ of order $k$ (we call this property a $s_k$-totally associativity). 
When $p=2$ and $q=1$, we obtain a $(2k+1)$-ary product on the space of bilinear maps on 
$E$ with values on $E$, which is identified to the cubic matrices. 
If we call a $l$-matrix a square tableau  with $l\times \cdots\times l$ entrances 
(if $l=3$ we have the cubic matrices and we speak about hypercubic matrices as soon as $l >3$), then
the $(2k+1)$-ary product on $T^p_q(E)$  gives a $(2k+1)$-product on the space of $(p+q)$-matrices. 
We describe also all 
these products which are $s_k$-totally associative. We compute the corresponding quadratic operads and their dual.

\end{abstract}

\medskip


\section{On $n$-ary associative algebras}

\subsection{Generalities: $n$-ary partially and totally associative algebras}

A $n$-ary algebra is a pair $(V,\mu)$ where $V$ is a vector space on a commutative field $\mathbb{K}$
of characteristic $0$ and $\mu$ a linear map 
\begin{equation*}
\mu :V^{\otimes n}\rightarrow V
\end{equation*}
where $V^{\otimes n}$ denotes the $n$-tensor product $V\otimes \cdots \otimes V$ ($n$ times).

A $n$-ary algebra is  {\it partially associative} if $\mu$ satisfies

\begin{equation}
\label{par}
\displaystyle{\sum_{p=0}^{n-1} (-1)^{p(n-1)} \mu \circ (I_p \otimes \mu \otimes I_{n-p-1})=0},
\end{equation}
for any $p=0, \cdots ,n-1$, where $I_j:V^{\otimes j}\rightarrow V^{\otimes j}$ is the identity map and $I_0\otimes \mu=\mu \otimes I_0=\mu.$
For example, if $n=2$ we have the classical notion of binary associative algebra.

A $n$-ary algebra is  {\it totally associative} if $\mu$ satisfies
\begin{equation}
\label{tot}
\displaystyle{ \mu \circ (\mu \otimes I_{n-1})
=  \mu \circ (I_p \otimes \mu \otimes I_{n-p-1})},
\end{equation}
for all $p=0, \cdots, n-1.$
For $n=2$, the two notions of partially and totally associativity coincide 
with the classical notion of associativity.
A totally associative $(2p)$-ary algebra is partially associative. A totally associative $(2p+1)$-ary
algebra is partially associative if and only if $\mu$ is $2$-step nilpotent (i.e. $\mu \circ_i \mu = 0$ for
any $i=1, \cdots, n$ with $\mu \circ_i \mu=\mu \circ (I_{i-1} \otimes \mu \otimes I_{n-i})$).

\noindent{\bf Remark. } Some applications of cubic or $n$-ary algebras in physic can be found in \cite{C.RT} or \cite{RT1}
and \cite{RT2}.

\subsection{Definition of $n$-ary $\sigma$-partially and $\sigma$-totally associative algebras}

 We can generalize Identities (\ref{par})  and 
(\ref{tot}) using actions of the symmetric group on $n$ elements 
$\Sigma_n.$ This generalization is in the spirite 
of the binary $\mathbb{K}[\Sigma_3]$-associative algebras
introduced and developped in \cite{G.R1} and \cite{G.R3}.

\begin{definition}
For a permutation $\sigma$ in $\Sigma_n$ define a linear map 
$$\Phi_\sigma^V: 
V^{\otimes^n} \rightarrow V^{\otimes^n}$$
by 
$$\Phi_\sigma^V(e_{i_1} \otimes \cdots \otimes e_{i_n})
=e_{i_{\sigma^{-1}(1)}} \otimes \cdots \otimes e_{i_{\sigma^{-1}(n)}}.$$
A $n$-ary algebra $(V,\mu)$ is 
{\bf $\sigma$-partially associative} if 
\begin{equation}
\label{sigmapar}
\displaystyle{\sum_{p=0}^{n-1} (-1)^{p(n-1)} 
(-1)^{p\varepsilon (\sigma)} 
\mu \circ (I_p \otimes( \mu\circ \Phi_{\sigma^p}^V)  
\otimes I_{n-p-1})=0},
\end{equation}
for all $p=0, \cdots, n-1,$

\noindent  and  {\bf $\sigma$-totally associative } if  
\begin{equation}
\label{sigmatot}
\displaystyle{ \mu \circ (\mu \otimes I_{n-1})= 
\mu \circ (I_p \otimes (\mu\circ \Phi_{\sigma^p}^V)  
 \otimes I_{n-p-1})},
\end{equation}
for all $p=0, \cdots, n-1.$
\end{definition}

\noindent {\bf Example} If $n=3$ and $\sigma=\tau_{12}$ is
the transposition exchanging $1$ and $2$ then a 
$\tau_{12}$-totally associative algebra satisfies
$$ \mu(\mu(e_1,e_2,e_3),e_4,e_5)=
\mu(e_1,\mu (e_3,e_2,e_4),e_5)=
\mu(e_1,e_2(\mu (e_3,e_4,e_5)),$$
and a $\tau_{12}$-partially associative algebra satisfies
$$ \mu(\mu(e_1,e_2,e_3),e_4,e_5)-
\mu(e_1,\mu (e_3,e_2,e_4),e_5)+
\mu(e_1,e_2(\mu (e_3,e_4,e_5))=0.$$

\section{A $(2p+1)$-ary product on the vector space of tensors $T_1^2(E)$}
 
\subsection{The tensor space $T_1^2(E)$}

Let $E$ be a finite dimensional vector space over a field $\mathbb{K}$ of characteristic $0$. We denote by
$T_1^2(E)=E \otimes E \otimes E^*$ the space of tensors covariant of degree $1$ and contravaviant of degree $2.$
The space $T_1^2(E)$ is identified to the space of linear maps 
$$\mathcal{L}(E \otimes E,E)=\left\{ \varphi: E \otimes E \rightarrow E \ \mbox{\rm linear} \right\}.$$  
Let $\left\{ e_1,\cdots ,e_n \right\}$ be a fixed basis of $E$. The structure constants $\left\{ C_{ij}^k \right\}$ 
of $ \varphi  \in T^2_1(E)$ are defined by
$$ \varphi(e_i \otimes e_j) = \displaystyle{\sum_{k=1}^n} C_{ij}^k e_k.$$

\begin{definition}
The dual map of $ \varphi \in T^2_1(E)$ is the tensor 
$\widetilde{ \varphi} \in T^1_2(E)\simeq \mathcal{L}(E,E\otimes E)   $ 
defined by 
$$\begin{array}{cccc}
\widetilde{\varphi}: & E & \rightarrow & E \otimes E \\
& e_k & \mapsto & \displaystyle{\sum_{1\leq i,j \leq n}} C_{ij}^k e_i \otimes e_j.
\end{array}$$
 \end{definition}
If $ \varphi$ is considered as a multiplication on $E$, then $\widetilde{\varphi}$ 
 is a coproduct. For example, if $\varphi$ is an associative product
then $\widetilde{\varphi}$ is the corresponding coassociative coproduct (often denoted by
$\Delta$). 

\subsection{Definition of a $3$-ary product on $T_1^2(E)$}

Let $\varphi_1,\varphi_2,\varphi_3$ be in $T_1^2(E).$ We define a $3$-ary 
product $\mu$ by
\begin{eqnarray}
\label{def}
\mu(\varphi_1,\varphi_2,\varphi_3 )=\varphi_1 \circ \widetilde{\varphi_2} \circ \varphi_3. 
\end{eqnarray}
As $\widetilde{\varphi_2}: E\rightarrow E \otimes E$, then 
$\varphi_1 \circ \widetilde{\varphi_2} \circ \varphi_3 \in T^2_1(E)$
and  $\mu$ is well defined. Let us compute its stucture constants. We denote by $C^k_{ij}(l)$ the 
structure constants of $\varphi_l$ $(l=1,2,3).$

$$\begin{array}{ll}
\mu (\varphi_1,\varphi_2,\varphi_3)(e_i \otimes e_j)=& 
\varphi_1 \circ \widetilde{\varphi_2} \circ \varphi_3(e_i \otimes e_j)\\
& = \displaystyle{\sum_{k=1}^n} C_{ij}^k(3) \varphi_1 \circ \widetilde{\varphi_2}
(e_k) \\
& = \displaystyle{\sum_{k=1}^n} \ \ \displaystyle{\sum_{1 \leq l,m \leq n}} 
C_{ij}^k(3) C_{lm}^k(2) \varphi_1 (e_l \otimes e_m) \\
& =\displaystyle{\sum_{t=1}^n} \  \displaystyle{\sum_{k=1}^n} \ \  \displaystyle{\sum_{1 \leq l,m \leq n}} 
C_{ij}^k(3) C_{lm}^k(2) C_{lm}^t(1) e_t.
\end{array}
$$ 
Thus if $\mu (\varphi_1,\varphi_2,\varphi_3)(e_i \otimes e_j)=
\displaystyle{\sum_{t=1}^n}A_{ij}^t(1,2,3)e_t$ we get 
$$A_{ij}^t(1,2,3)=\displaystyle{\sum_{1\leq k,l,m \leq n}} C_{ij}^k(3) C_{lm}^k(2) C_{lm}^t(1).$$

\begin{proposition}
The $3$-ary product in $T_1^2(E)$ given by
$$\mu ( \varphi_1,\varphi_2,\varphi_3)=\varphi_1 \circ \widetilde{\varphi_2} \circ \varphi_3$$
satisfies
$$\begin{array}{ll}
\mu (\mu ( \varphi_1,\varphi_2,\varphi_3), \varphi_4,\varphi_5) & =
\mu (\varphi_1, \mu ( \varphi_4,\varphi_3,\varphi_2),\varphi_5) \\
& =\mu ( \varphi_1,\varphi_2,\mu (\varphi_3, \varphi_4,\varphi_5)),
\end{array}$$
that is this product is $\tau_{13}$-totally associative.
\end{proposition}

\noindent {\it Proof.} We have
$$\begin{array}{ll}
\mu ( \mu (\varphi_1,\varphi_2,\varphi_3),\varphi_4,\varphi_5)
(e_i \otimes e_j) & 
=(\varphi_1 \circ \widetilde{\varphi_2} \circ \varphi_3) \circ \widetilde{\varphi_4} \circ 
\varphi_5 (e_i \otimes e_j) \\
& = \displaystyle{\sum_t} \left[ \displaystyle{\sum_{k,l,m}} C_{ij}^k(5) C_{lm}^k(4) A_{lm}^t(1,2,3)\right]e_t \\
& = \displaystyle{\sum_t} \left[ \displaystyle{\sum_{k,l,m}} \, \displaystyle{\sum_{u,r,s}} \left(
C_{ij}^k(5) C_{lm}^k(4) C_{lm}^u(3) 
C_{rs}^u(2) C_{rs}^t(1) \right) \right]e_t.
\end{array}
$$ 
Thus the structure constant $A_{ij}^t((1,2,3),4,5)$ of this tensor is
$$A_{ij}^t((1,2,3),4,5)=\displaystyle{\sum_{1 \leq
\begin{array}{l}
k,l,m \\
u,r,s
\end{array}
\leq n }} C_{ij}^k(5) C_{lm}^k(4) C_{lm}^u(3) 
C_{rs}^u(2) C_{rs}^t(1).
$$
Similary
$$\begin{array}{ll}
\mu ( \varphi_1,\varphi_2,\mu(\varphi_3,\varphi_4,\varphi_5))
(e_i \otimes e_j) & =
\varphi_1 \circ \widetilde{\varphi_2} \circ ( \varphi_3 \circ \widetilde{\varphi_4} \circ 
\varphi_5 )(e_i \otimes e_j) \\
& = \displaystyle{\sum_t} \left[ \displaystyle{\sum_{u,r,s}} A_{ij}^u(3,4,5)C_{rs}^u(2) C_{rs}^t(1)\right]e_t \\
& = \displaystyle{\sum_t} \left[ \displaystyle{\sum_{u,r,s}} \left( \displaystyle{\sum_{k,l,m}}
C_{ij}^k(5) C_{lm}^k(4) C_{lm}^u(3) \right)
C_{rs}^u(2) C_{rs}^t(1) \right]e_t.
\end{array}
$$ 
Thus
$$A_{ij}^t (1,2,(3,4,5))=\displaystyle{
\sum_{
\begin{array}{l}
k,l,m \\
u,r,s
\end{array}
}} C_{ij}^k(5) C_{lm}^k(4) C_{lm}^u(3) C_{rs}^u(2) C_{rs}^t(1),
$$
and 
$$A_{ij}^t (1,2,(3,4,5))=A_{ij}^t ((1,2,3),4,5).$$
We also have 
$$\begin{array}{ll}
\mu ( \varphi_1,\mu (\varphi_2,\varphi_3,\varphi_4),\varphi_5)
(e_i \otimes e_j) & =\varphi_1 \circ (\widetilde{\varphi_2 \circ  \widetilde{\varphi_3} \circ \varphi_4}) \circ 
\varphi_5 (e_i \otimes e_j) \\
& = \displaystyle{\sum_t} \left[ \displaystyle{\sum_{k,l,m}} C_{ij}^k(5) A_{lm}^k(2,3,4) C_{lm}^t(1)\right]e_t \\
& = \displaystyle{\sum_t} \left[ \displaystyle{\sum_{k,l,m}}  \displaystyle{\sum_{u,r,s}} 
C_{ij}^k(5) C_{lm}^u(4) C_{rs}^u(3) 
C_{rs}^k(2) C_{lm}^t(1)  \right]e_t,
\end{array}
$$ 
and
$$A_{ij}^t (1,(2,3,4),5)=\displaystyle{
\sum_{
\begin{array}{l}
k,l,m \\
u,r,s
\end{array}
}} C_{ij}^k(5) C_{lm}^u(4) C_{rs}^u(3) C_{rs}^k(2) C_{lm}^t(1).$$
This shows that
$$A_{ij}^t ((1,2,3),4,5)=A_{ij}^t (1,(4,3,2),5).$$

\noindent {\bf Remarks.} 

\noindent 1. We can define in this way other non equivalent products by:
$$\left\{
\begin{array}{l}
\mu_2 ( \varphi_1,\varphi_2,\varphi_3)=  \varphi_3 \circ \widetilde{\varphi_2} \circ \varphi_1, \\
\mu_3 ( \varphi_1,\varphi_2,\varphi_3)=  \varphi_1 \circ \widetilde{\varphi_2} \circ \, ^t\varphi_3, \\
\mu_4 ( \varphi_1,\varphi_2,\varphi_3)=  \varphi_3 \circ \widetilde{\varphi_2} \circ \, ^t\varphi_1,
\end{array}
\right.
$$
where  $ ^t \varphi (e_i \otimes e_j)=\varphi (e_j \otimes e_i).$
 
\noindent 2. If we identify a tensor $\varphi$ to its structure constants $\left\{ C_{ij}^k \right\}$
and if we consider the family $\left\{ C_{ij}^k \right\}$ as a cubic matrix $\left\{ C_{ijk} \right\}$
with $3$-entries, the product $\mu$ on $T^2_1(E)$ gives a $3$-ary product on the cubic matrices. 
This last product has been studied in \cite{A.K.L.S}.

\subsection{A $(2k+1)$-ary product on $T^2_1(E)$}

Let $\varphi_1 , \cdots ,\varphi_{2k+1}$ be in $T^2_1(E)$. We define a $(2k+1)$-ary product 
$\mu_{2k+1}$ on $T^2_1(E)$ putting
$$\mu_{2k+1}(\varphi_{1} , \cdots, \varphi_{2k+1} )= \varphi_{1} \circ \widetilde{\varphi_{2}} \circ
\cdots \circ \varphi_{2k-1} \circ \widetilde{\varphi_{2k}} \circ \varphi_{2k+1} .$$
Let $s_k$ be the permutation of $\Sigma_{2k+1}$ defined by 
$$s_k(1,2,\cdots , 2k+1)=(2k+1,2k, \cdots , 2,1),$$
that is $s_k= \tau_{1 \, 2k+1} \circ \tau_{2 \, 2k} \circ \cdots \circ \tau_{k-1 \, k+1}
=\Pi_{i=i}^k \tau_{i\, 2k+1-i}$. It satisfies 
$(s_k)^{2p}=Id$ and  $(s_k)^{2p+1}=s_k$ for any $p$ (it is a symmetry).

Recall that the $(2k+1)$-ary product $\mu_{2k+1}$ is a $s_k$-totally associative product if
$$\begin{array}{l}
\mu_{2k+1} \circ (\mu_{2k+1} \otimes I_{2k} )  
= \mu_{2k+1} \circ (I_p \otimes (\mu_{2k+1} \circ \Phi_{s_k^p})\otimes I_{2k-p}), 
\end{array}$$
for $p=1,\cdots,2k.$
In particular, we have
$$\mu_{2k+1} \circ (\mu_{2k+1} \otimes I_{2k} )  = 
 \mu_{2k+1} \circ (I_{2q} \otimes \mu_{2k+1} \otimes I_{2k-2q}),$$
for any $q=1, \cdots ,k.$

\begin{proposition} The product $\mu_{2k+1}$ is  $s_k$-totally associative.
\label{propo}
\end{proposition}

\noindent {\it Proof.} In fact if we put
$$\mu_{2k+1}(\varphi_1 , \cdots ,\varphi_{2k+1})(e_i \otimes e_j)= \displaystyle{\sum_t }A_{ij}^t(1,2,\cdots,2k+1)e_t,$$
then $A_{ij}^t(1,2,\cdots,2k+1)=$
$$\displaystyle{\sum_
{\begin{array}{c}
a_1,\cdots,a_{k+1}\\
k_1,\cdots ,k_{k}
\end{array} }} C_{ij}^{k_1}(2k+1) C_{a_1a_2}^{k_1}(2k) C_{a_1a_2}^{k_2}(2k-1) 
\cdots C_{a_{2k-1}a_{2k}}^{k_k}(2) C_{a_{2k-1}a_{2k}}^t(1). 
$$
More precisely the line of superscripts is 
$$(k_1,k_1,k_2,k_2, \cdots , k_k, k_k, t),$$
and the line of subscripts
$$((i,j), (a_1,a_2),(a_1,a_2),(a_3,a_4),(a_3,a_4), \cdots , (a_{2k-1},a_{2k}),(a_{2k-1},a_{2k})).$$
Let us consider 
$$ \mu_{2k+1} \circ (I_{l} \otimes (\mu_{2k+1} \circ \Phi_{s_k^l})\otimes I_{2k-l})
(\varphi_1 , \cdots ,\varphi_{4k+1})(e_i \otimes e_j)= 
\sum B_{ij}^t e_t.$$ 
Thus for $l=2r$, we get
$$
\begin{array}{ll}
\medskip
B_{ij}^t=&\sum C_{ij}^{k_1}(4k+1) C_{a_1a_2}^{k_1}(4k) C_{a_1a_2}^{k_2}(4k-1) 
\cdots C_{a_{2k-2r-1}a_{2k-2r}}^{k_{k-r}}(2k+l+2)\\
&A_{a_{2k-2r-1}a_{2k-2r}}^{k_{k-r+1}}(l+1,\cdots,2k+l+1)C_{a_{2k-2r+1}a_{2k-2r+2}}^{k_{k-r+1}}(l)\cdots
 C_{a_{4k-1}a_{4k}}^t(1),
\end{array} $$
such that the line of superscripts is 
$$(k_1,k_1,k_2,k_2, \cdots ,k_{k-r},h_1,h_1,\cdots,h_k,h_k,k_{k-r+1},k_{k-r+1},\cdots, k_k, k_k, t),$$
where the terms $h_1,\cdots,h_k,k_{k-r+1}$ correspond to the factor $A_{a_{2k-2r-1}a_{2k-2r}}^{k_{k-r+1}}(l+1,\cdots,2k+l+1)$.
Such a line is the same as the line of superscripts of 
$$\mu_{2k+1} \circ  (\mu_{2k+1} \otimes I_{2k})
(\varphi_1 , \cdots ,\varphi_{4k+1})(e_i \otimes e_j).$$ 
The line of subscripts is
$$
\begin{array}{l}
((i,j), (a_1,a_2),(a_1,a_2),\cdots , (a_{2k-2r-1},a_{2k-2r}),(a_{2k-2r-1},a_{2k-2r}) ,(\beta_1 \beta_2),\cdots,
(\beta _{2k-1},\beta _{2k}),\\(a_{2k-2r-1},a_{2k-2r}),\cdots, (a_{4k-1},a_{4k})).
\end{array}
$$
So
$$\mu_{2k+1} \circ ( (\mu_{2k+1} \otimes I_{2k})=\mu_{2k+1} \circ (I_{l} \otimes (\mu_{2k+1} \circ \Phi_{s_k^l})
\otimes I_{2k-l}),$$
for $l=2r.$
Assume now that $l=2r+1$. In this case $B_{ij}^t$ is of the form
$$\displaystyle{\sum} \cdots C_{a_{2k-2r-1}a_{2k-2r}}^{k_{k-r+1}}(2k+l+2)
A_{a_{2k-2r+1}a_{2k-2r+2}}^{k_{k-r+1}}(2k+l+1,\cdots,l+1)
C_{a_{2k-2r+1}a_{2k-2r+2}}^{k_{k-r+1}}(l)\cdots.$$
We find the same list of exponents and of indices that for
$\mu_{2k+1} \circ  (\mu_{2k+1} \otimes I_{2k})$. This finishes the proof.

\medskip

\noindent{\bf Consequences.}

1. The product $\mu_{2k+1}$ on $T^2_1(E)$ induces directly a $(2k+1)$-ary products on cubic matrices.

2. All the other products which are $s_k$-totally associative corresponds to
$$\left\{
\begin{array}{lll}
\mu_{2k+1}^2 (\varphi_1 , \cdots ,\varphi_{2k+1}) & = & \varphi_{2k+1} \circ \widetilde{\varphi_{2k}}  \circ  \cdots 
\widetilde{\varphi_{2}} \circ \varphi_{1}, \\
\mu_{2k+1}^3 (\varphi_1 , \cdots, \varphi_{2k+1}) & = & \mu_{2k+1} ( \, ^t\varphi_1 ,\varphi_{2}, \cdots, \varphi_{2k+1}), \\ 
\mu_{2k+1}^4 (\varphi_1 , \cdots ,\varphi_{2k+1}) & = & \mu_{2k+1}^2 ( \varphi_1 , \cdots, \varphi_{2k} ,\, ^t\varphi_{2k+1}).
\end{array}
\right.
$$ 
and more generally
$$\mu_{2k+1} ( \, ^t\varphi_1 , \varphi_2, ^t\varphi_3,\cdots, \varphi_{2k+1})$$ 
or 
$$\mu_{2k+1} ( \, ^t\varphi_1 , \varphi_2, ^t\varphi_3,\cdots, \varphi_{2k+1}).$$

\section{Generalisation: a $(2k+1)$-ary product on $T^p_q(E)$}

\subsection{The vector space $T_q^p(E)$}

Let $E$ be a finite $m$-dimensional $\mathbb{K}$-vector space. 
 The vector space $T_q^p(E)$ is the space of tensors which are contravariant of degree $p$ and  covariant of degree $q.$  
In  $\left\{ e_1, \cdots ,e_m \right\}$ is a fixed basis of $E$, a tensor $t$ of 
$T_q^p(E)$ is written
$$ t= \displaystyle{
 \sum_{\begin{array}{c}
1 \leq i_k, j_l \leq n \\
1 \leq k \leq p \\
1 \leq l \leq q 
\end{array}}
}
t_{i_1, \cdots , i_p}^{j_1, \cdots , j_q} e_{i_1}\otimes \cdots \otimes e_{i_p}\otimes e^{j_1}\otimes \cdots 
\otimes  e^{j_q}$$
where $(e^1, \cdots, e^m)$ is the dual basis of $(e_1, \cdots , e_m).$
As
$$ T_q^p(E)=T_0^p(E) \otimes T^0_q(E),$$ then the tensor space 
$$T(E)=\displaystyle{\sum_{p,q=0}^\infty }T_q^p(E)$$
is an associative algebra with product 
$$\begin{array}{ccc}
T_q^p(E) \times T_{m}^{l}(E) & \rightarrow  & T_{q+m}^{p+l}(E) \\
(K,L) & \mapsto & K \otimes L
\end{array}.$$
But this product is not internal on each component $T_q^p(E).$ In this section we will define  internal 
$(2p-1)$-ary-product on the components. 

The vector space $T^p_q(E)$ is isomorphic to the space $\mathcal{L}(E ^{\otimes p},E^{\otimes q})$ of linear maps 
$$t: E^{\otimes^p} \rightarrow E^{\otimes^q} .$$

We define the structure constants by 
$$t(e_{i_1} \otimes \cdots \otimes e_{i_p})=
\sum C_{i_1 \, \cdots \, i_p}^{j_1 \, \cdots \, j_q} e_{j_1} \otimes \cdots \otimes e_{j_q}.$$
For such a map we define $\widetilde{t}$ by 
$$\begin{array}{cccc}
\widetilde{t}: &  E^{\otimes^q} & \rightarrow  & E^{\otimes^p} \\
& (e_{j_1} \otimes \cdots \otimes e_{j_q}) & \mapsto  & 
\sum C^{j_1 \, \cdots \, j_q}_{i_1 \, \cdots \, i_p} e_{i_1} \otimes \cdots \otimes e_{i_p}.
\end{array}
$$

\subsection{A $(2k+1)$-ary product on $T^p_q(E)$}

\begin{definition} The map $\mu$ defined by:
\begin{eqnarray}
\label{pdt}
\mu(\varphi_1 , \cdots ,\varphi_{2k+1})=
\varphi_{2k+1} \circ \widetilde{\varphi_{2k}} \circ \varphi_{2k-1} \circ  \cdots \circ
\widetilde{\varphi_{2}} \circ \varphi_1 , 
\end{eqnarray}
 for any $\varphi_1, \cdots ,\varphi_{2k+1} \in T^s_r(E)$ 
is a  $(2k+1)$-ary product on  $T^s_r(E)$.
\end{definition}

We take an odd number of map $\varphi_i$  so we get
compostions of 
$\widetilde{\varphi_{j+1}} \circ \varphi_j: E^{\otimes^p} \rightarrow E^{\otimes^p}$  for $j=1,\cdots, 2k-1$ 
and finally compose with $\varphi_{2k-1}: E^{\otimes^p} \rightarrow E^{\otimes^q}$ so $\mu$ is well defined.

\begin{proposition}
The $(2k+1)$-ary product $\mu$ on  $T^p_q(E)$ defined by (\ref{pdt}) is  $s_k$-totally associative.
\end{proposition}

\noindent {\it Proof.} The proof is similar to the proof of Proposition \ref{propo} 
concerning an $(2k+1)$-ary product on $T^2_1(E)$.  
In fact we have 
$$\mu (\varphi_1, \cdots \varphi_{2p+1})
(e_{i_1} \otimes \cdots \otimes e_{i_p})=\sum
A_{i_1 \cdots i_p}^{r_1 \cdots r_q} e_{r_1} \otimes \cdots \otimes e_{r_q}, $$ 
and $$A_{i_1 \cdots i_p}^{r_1 \cdots r_q}=C_{i_1 \cdots i_p}^{j_1 \cdots j_q}(2k+1)C_{l_1 \cdots l_p}^{j_1 \cdots j_q}(2k)
C_{l_1 \cdots l_p}^{m_1 \cdots m_q}(2k-1) \cdots C_{s_1 \cdots s_p}^{r_1 \cdots r_q}(1),$$
that is the line of superscripts  is
$$ (j_1 \cdots j_q)(j_1 \cdots j_q)(m_1 \cdots m_q)(m_1 \cdots m_q)  \cdots (n_1 \cdots n_q) 
(n_1 \cdots n_q) (r_1 \cdots r_q),$$
and the line of subscripts is
$$ (i_1 \cdots i_p)(l_1 \cdots l_p)(l_1 \cdots l_p)  \cdots  (s_1 \cdots s_p) (s_1 \cdots s_p). $$
Using the same arguments that in Proposition \ref{propo}, changing pairs by $p$-uples and $q$-uples, we obtain 
the announced result.

\medskip

\noindent {\bf Remark.} We can also use the same trick 
that in Consequences 2. to find others $s_k$-totally associative products on $T^p_q(E).$   

\bigskip

\noindent{\bf Applications.  } This product can be translated as a product of "hypercubic matrices" that is 
square tableau of length
$p+q$. This generalizes in a natural way the classical associative product of matrices. 

\section{Current $(2k+1)$-ary $s_k$-totally associative algebras}

The name refers to current Lie algebras which are Lie algebras of the form $L \otimes A$ where $L$ is a
Lie algebra and $A$ is a associative commutative algebra, equipped with bracket
$$[x \otimes a, y \otimes b]_{L \otimes A } = [x, y]_{L} \otimes ab.$$ We want to generalize this notion to 
$(2k+1)$-ary $s_k $-totally associative algebras. The problem is to find a category of $(2k+1)$-ary algebras such that 
its tensor product with a  $(2k+1)$-ary $s_k$-totally associative 
algebra gives a $(2k+1)$-ary $s_k$-totally associative algebra with obvious operation on the tensor product. 
Such a tensor product will be called
{\it  current $(2k+1)$-ary $s_k$-totally associative algebra}. We first focus on the 
ternary case and $s_1=\tau_{13}$.
 
Let $(V,\mu)$ be a $3$-ary algebra where $\mu$ is a $\tau_{13}$-totally associative product 
on $V$ (for example $V=T_1^2(E)$ and $\mu$ is defined by (\ref{def}) ) so $\mu$ satisfies Equations (\ref{sigmatot})
for $\sigma=\tau_{13}$, that is,  
$$ \mu(\mu(e_1,e_2,e_3),e_4,e_5)=
\mu(e_1,\mu (e_4,e_3,e_2),e_5)=
\mu(e_1,e_2(\mu (e_3,e_4,e_5)),$$
for any $e_1,e_2,e_3$ in $V.$
Let $(W,\tilde{\mu})$ be a $3$-ary algebra. Then the tensor algebra 
$\left( V\otimes W,\mu\otimes \tilde{\mu}\right)$
is a $3$-ary $\tau_{13}$-totally associative algebra if and only if
$$(\mu\otimes \tilde{\mu})(v_1 \otimes w_1 \otimes v_2 \otimes w_2 \otimes v_3 \otimes w_3)= \mu(v_1,v_2,v_3) \otimes 
\tilde{\mu}(w_1,w_2,w_3)$$
satisfies the $\tau_{13}$-totally associativity relation. But
$$\left\{
\begin{array}{l}
 (\mu\otimes \tilde{\mu})\circ (\mu\otimes \tilde{\mu} \otimes I_4) =  \mu \circ (\mu \otimes I_2) \otimes \tilde{\mu }
\circ  (\tilde{\mu} \otimes I_2), \\
 (\mu\otimes \tilde{\mu})\circ ( I_2 \otimes (\mu\otimes \tilde{\mu}) \circ 
\Phi_{\tau_{13}}^{V \otimes W} \otimes I_2)
 =  \mu \circ (I \otimes \mu \circ \Phi_{\tau_{13}}^{V } \otimes I) \otimes \tilde{\mu } \circ (I \otimes \tilde{\mu} 
\circ \Phi_{\tau_{13}}^W \otimes I),\\ 
 (\mu\otimes \tilde{\mu})\circ (I_4  \otimes \mu\otimes \tilde{\mu}) =  \mu \circ (I_2 \otimes \mu) 
\otimes \tilde{\mu } \circ  (I_2 \otimes \tilde{\mu}), 
\end{array}
\right.$$
then $(\mu \otimes \tilde{\mu}) \circ (\mu \otimes \tilde{\mu} \otimes I_4)-(\mu \otimes \tilde{\mu}) \circ
(I_4 \otimes \mu \otimes \tilde{\mu})=0$
is equivalent to
\begin{eqnarray}
\label{equa1}
\mu \circ (\mu \otimes I_2) \otimes \tilde{\mu} \circ (\tilde{\mu} \otimes I_2)-\mu \circ (I_2 \otimes \mu) \otimes 
\tilde{\mu} \circ (I_2 \otimes \tilde{\mu})=0.
\end{eqnarray}
But $\mu \circ (\mu \otimes I_2)=\mu \circ (I_2 \otimes \mu)$. Thus Equation (\ref{equa1}) is equivalent to
$$\mu \circ (\mu \otimes I_2) \otimes 
\left[ \tilde{\mu} \circ (\tilde{\mu} \otimes I_2) - \tilde{\mu} \circ(I_2 \otimes \tilde{\mu}) \right] =0,$$
and $$\tilde{\mu} \circ (\tilde{\mu} \otimes I_2)= \tilde{\mu} \circ (I_2 \otimes \tilde{\mu}).$$
Similary
$$\begin{array}{l}
(\mu\otimes \tilde{\mu})\circ (\mu\otimes \tilde{\mu} \otimes I_4)-
(\mu\otimes \tilde{\mu})\circ (I_4  \otimes \mu\otimes \tilde{\mu})\\
\mu \circ (\mu \otimes I_2) \otimes 
\left[ \tilde{\mu} \circ (\tilde{\mu} \otimes I_2) - \tilde{\mu} 
\circ(I \otimes \tilde{\mu}\circ \Phi^W_{\tau_{13}} \otimes I) \right] =0,
\end{array}$$
which leads to 
$$\tilde{\mu} \circ (\tilde{\mu} \otimes I_2) = \tilde{\mu} 
\circ(I \otimes \tilde{\mu}\circ \Phi^W_{\tau_{13}} \otimes I).$$
So $\mu \otimes \tilde{\mu}$ is $\tau_{13}$-totally associative if and only if 
$\tilde{\mu}$ is $\tau_{13}$-totally associative.

\begin{proposition}
\label{3tensorproduct}
Let $(V,\mu)$ be a $3$-ary $\tau_{13}$-totally associative algebra and $(W,\tilde{\mu})$ be a $3$-ary algebra.
Then $(V\otimes W,\mu \otimes \tilde{\mu})$ is a  $3$-ary $\tau_{13}$-totally associative algebra
if and only if $(W,\tilde{\mu})$ is also of this type.
\end{proposition}

This result can be extended for $(2k+1)$-ary $s_k$-totally associative algebras.

\begin{proposition}
Let $(V,\mu)$ be a $(2k+1)$-ary $s_k$-totally associative algebra and $(W,\widetilde{\mu})$ be a $(2k+1)$-ary algebra.
Then $(V\otimes W,\mu \otimes \tilde{\mu})$ is a  $(2k+1)$-ary $s_k$-totally associative algebra
if and only if $(W,\tilde{\mu})$ is also of this type.
\end{proposition}

\noindent {\it Proof.} The product $\mu$ is a $(2k+1)$-ary $s_k$-totally associative product so satisfies 
$$\begin{array}{lll}
\mu \circ (\mu \otimes I_{2k} ) & = &
 \mu \circ (I_{2q} \otimes \mu \otimes I_{2k-2q})\\
 & = &  \mu \circ (I_{2q+1} \otimes \mu \circ \Phi_{s_k^q}^V \otimes I_{2k-2q-1}),
\end{array}$$
for any $q=0, \cdots ,k.$
The system
$$\begin{array}{l}
(\mu \otimes \widetilde{\mu}) \circ ((\mu \otimes \widetilde{\mu}) \otimes I_{4k} )  - 
 (\mu \otimes \widetilde{\mu}) \circ (I_{4q} \otimes (\mu \otimes \widetilde{\mu})
\circ \Phi_{s_k^q}^{V \otimes W} 
\otimes I_{4k-2q})=\\
\mu \circ (\mu \otimes I_{2k} ) 
\otimes \widetilde{\mu} \circ (\widetilde{\mu} \otimes I_{2k}) - 
 \mu \circ (I_{q} \otimes \mu \circ \Phi_{s_k^q}^{V } \otimes I_{2k-q})
\otimes  \tilde{\mu} \circ (I_{q} \otimes \widetilde{\mu} \circ \Phi_{s_k^q}^{W } \otimes I_{2k-q})=0,
\end{array}$$
for any $q=0, \cdots ,k$ is equivalent to
$$\begin{array}{l}
\mu \circ (\mu \otimes I_{2k} ) 
\otimes \left[ \widetilde{\mu} \circ (\widetilde{\mu} \otimes I_{2k}) - 
 \tilde{\mu} \circ (I_{q} \otimes \widetilde{\mu} \circ \Phi_{s_k^q}^{W } \otimes I_{2k-q}) \right] =0,
\end{array}$$
for any $q=0, \cdots ,k.$ Then $\mu \otimes \tilde{\mu}$ is  $(2k+1)$-ary $s_k$-totally associative if and only if
$$\begin{array}{l}
\widetilde{\mu} \circ (\widetilde{\mu} \otimes I_{2k}) - 
 \tilde{\mu} \circ (I_{q} \otimes \widetilde{\mu} \circ \Phi_{s_k^q}^{W } \otimes I_{2k-q}) =0\\
\end{array}$$
for any $q=0, \cdots ,k$ that is  $\tilde{\mu}$ is  a $(2k+1)$-ary $s_k$-totally associative product.

\section{The operads $\pa$, $\tota$ }

\subsection{On the operad $\pa$}

We denote by $\pa$ the quadratic operad of 
$3$-ary -i.e. ternary- partially associative algebras (with operation in degree $0$). 
In \cite{G.R} we compute the free $3$-ary partially associative algebra 
based on a finite dimensional vector space $V.$ In  \cite {R} we notice that the dual 
operad is in general defined in the graded  framework, compute it, as the knowledge of the dual is fundamental
to study if the operad is Kozsul or not. We prove in 
 \cite{R} that $\pa$ is not Koszul.
Note that this result contradicts some affirmations of the Koszulity of the operad $\pa$. 
This confusion can be explained by observing the general case of the operad $\npa$ for $n$-ary 
partially associative algebras with operation of degree $0$. 
If $n$ is even (\cite{Gne}),   $\npa$ is Koszul  and the dual operad  $\npa^!$ is the operad 
$\nta$ for  $n$-ary 
totally associative algebras with operation of degree $0$ (which is also Koszul).
But if $n=2k+1$, the operad   $\npa^!$ is not  $\nta$ but  $\ntad$ 
for totally associative algebras  with {\it operation of degree $1$} and this operad is not Koszul 
(see\cite{R}). As a consequence we deduce that for $n$ odd, the operadic cohomology (which always exits)
is not the cohomology which governs deformations (which also always exits, contrary to what is written in \cite{GoM}). Remark that in \cite{G.R} we have also defined a cohomology of Hochschild 
type for $3$-ary partially associative algebras with some extra conditions.

\subsection{The operad $\tota$}

We denote by $\tota$ the quadratic operad for $3$-ary $\tau_{13}$-totally 
associative algebras that is satisfying 
Relation (\ref{sigmapar}) for $\sigma=\tau_{13}$. 
Let $\mu$ be a $3$-ary multiplication,  and 
$$ E_{\tota}(m)=\left\{
\begin{array}{ll}
<\mu >\simeq \mathbb{K}[\Sigma_3], & \mbox{\rm if $m=3$ and} \\
0, & \mbox{\rm if $m\neq 3.$ }
\end{array}
\right. $$
We simply say that $E_{\tota}=E_{\tota}(3).$
The ideal of relation is generated by the $\mathbb{K}[\Sigma_{5}]$-closure  $R_{\tota}$ of the 
$\tau_{13}$-associativity
$$\left\{
\begin{array}{l}
r_1=\mu (\mu \otimes I_2)-\mu (I \otimes \mu \cdot \tau_{13} \otimes I), \\
r_2=\mu (\mu \otimes I_2)-\mu (I_2 \otimes \mu ) ,
\end{array}
\right.$$
where $\mu \cdot \sigma =\mu \circ \Phi_\sigma$ for $\sigma \in \Sigma_3$.

If $\Gamma(E_{\tota})$ denotes the free operad generated by $E_{\tota},$ we get that
$R_{\tota}\subset \Gamma(E_{\tota})(5).$
The operad for $3$-ary $\tau_{13}$-totally associative algebras is then the quadratic $3$-ary operad
$$\tota= \Gamma(E_{\tota})/(R_{\tota}),$$
that is $\tota(m)= \Gamma(E_{\tota})(m)/(R_{\tota})(m).$

\subsection{The current operad $\widetilde{\tota}$}
 In \cite{G.R4} we have defined, for a quadratic  operad $\mathcal{P}$, the current operad 
$\widetilde{\mathcal{P}}$ that is, the maximal operad $\widetilde{\mathcal{P}}$ such that the tensor product of a 
$\mathcal{P}$-algebra $A$ and a $\widetilde{\mathcal{P}}$-algebra $B$ is a $\mathcal{P}$-algebra with
the usual product on $A \otimes B$.
Let us compute  $\widetilde{\tota}.$

\begin{proposition}The current operad of the operad $\tota$ 
is $\tota$ itself that is
$$\widetilde{\tota}=\tota.$$ 
\end{proposition}
\pp This result follows from the Proposition \ref{3tensorproduct}.

\subsection{The dual operad $\tota^!$ }
For $n$-ary quadratic operad $\mathcal{P}=\Gamma(E)/(R)$ with $E=E(n)$, the dual (quadratic $n$-ary)
operad is defined as follows
$$\mathcal{P}^!= \Gamma(\overline{E})/(R^{\perp}),$$
where $\overline{E}=\uparrow^{n-2} E^\# \otimes sgn_n,$ $\uparrow^{n-2}$ denotes the suspension iterated $(n-2)$ times,
$\#$ the linear dual and $R^{\perp}\subset \Gamma(\overline{E})(2n-1)$ is the annihilator of 
$R\subset \Gamma(E)(2n-1)$ 
with respect to the pairing between $\Gamma(\overline{E})(2n-1)$  and $\Gamma(E)(2n-1).$ 

\begin{proposition} The dual operad of $\tota$ is   
$$\tota^!=\pad,$$ that is the operad for $\tau_{13} $-partially associative algebras with operation 
in degree $1$. 
\end{proposition}

\noindent {\it Proof.} The operad $\mathcal{P}=\tota$ is the quadratic operad defined by 
$$\mathcal{P}=\Gamma (E) / (R),$$ where $\mu$ a ternary operation of degree $0,$ $\Gamma (E)$
the free operad generated by $E=<\mu>$ and  $R \subset \Gamma(E)$ is the    generated as $\mathbb{K}[\Sigma_5]$-module
by the relations
$$\left\{
\begin{array}{l}
\mu(\mu \otimes I_2)-\mu(I \otimes \mu \cdot \tau_{13} \otimes I), \\
 \mu(\mu \otimes I_2)-\mu(I_2 \otimes \mu ).
\end{array} 
\right.$$ 
We consider 
$$\mu \circ_{s} \mu= \mu (I_{s-1} \otimes \mu \otimes I_{3-s}), $$
which "plugs" $\mu$ into the $s$-st input of $\mu$ and 
$$ (f \cdot \sigma) (i_1,i_2, \cdots , i_m)=f(  i_{\sigma^{-1}(1)},i_{\sigma^{-1}(2)},\cdots ,i_{\sigma^{-1}(m)}),$$
if $f \in \Gamma(\mu)(m), \sigma \in \Sigma_m$.

We get $\overline{E}(3)=<\alpha>$ where $\alpha$ is a ternary operation of degree $1$ satisfying 
$< \mu, \alpha >=1 .$  The pairing between $\Gamma(E)(5)$ and $\Gamma(\overline{E})(5)$ is given by
$$\begin{array}{l}
<(\mu \circ_j \mu )(i_1,i_2,i_3,i_4,i_5),(\alpha \circ_j \alpha>(i_1,i_2,i_3,i_4,i_5) ) \\
\qquad \qquad  = <\mu,\alpha> sgn_5 
\left(
\begin{array}{ccccc}
1 & 2 & 3 & 4 & 5 \\
i_1 & i_2 & i_3 & i_4 & i_5
\end{array}
\right) 
  =
sgn_5 
\left(
\begin{array}{ccccc}
1 & 2 & 3 & 4 & 5 \\
i_1 & i_2 & i_3 & i_4 & i_5
\end{array}
\right) ,
\end{array}
$$  
for $j=1,2,3.$ So
$$\begin{array}{l}
<(\mu \circ_1 \mu -\mu \circ_2 \mu  \cdot \tau_{13} )(1,2,3,4,5),
(\alpha\circ_1 \alpha -\alpha \circ_2 \alpha \cdot \tau_{13} +\alpha \circ_3 \alpha)  (1,2,3,4,5) > \\
\qquad = <\mu \circ_1 \mu,\alpha \circ_1 \alpha> +  <\mu \circ_2 \mu \cdot \tau_{13},\alpha \circ_2 \alpha \cdot \tau_{13}> \\
\qquad = 1-<\mu \circ_2 \mu ,\alpha \circ_2 \alpha>=1-1=0, \\
 \\
<(\mu \circ_1 \mu -\mu \circ_3 \mu)  (1,2,3,4,5),
(\alpha\circ_1 \alpha -\alpha \circ_2 \alpha \cdot \tau_{13} +\alpha \circ_3 \alpha)  (1,2,3,4,5) > \\
\qquad = <\mu \circ_1 \mu,\alpha \circ_1 \alpha> -  <\mu \circ_3 \mu ,\alpha \circ_3 \alpha >
=1-1=0.
\end{array}
$$ 
The dual operad is then the quadratic operad
$$\mathcal{P}^!=\Gamma (\alpha)/(R^{\perp}),$$ with $\alpha$ ternary operation of degree $1$ and 
$R^{\perp}$ generated by 
$$\begin{array}{l}
\alpha(\alpha \otimes I_2)-\alpha(I \otimes \alpha \cdot \tau_{13} \otimes I)
+\alpha(I_2 \otimes \alpha ).
\end{array}$$ So this operad is the operad of ternary $\tau_{13}$-partially associative algebras
with operations of degree $1.$ 
 
\medskip

\noindent {\bf Remark.} A direct computation similar to \cite{R} shows that 
$$dim \mathcal{P}(3)=6, \, dim \mathcal{P}(5)= 5! \, , \,  dim \mathcal{P}(7)=7! \quad .$$
The generating function of $\mathcal{P}$ is similar to the generating function of $\ta.$ 
Likewise the generating function of $\pad$ is the generating function of $\pd.$
>From \cite{H} the operads $\ta$ and $\pd$ are Koszul. We conclude that $\tota$ is Koszul.


\bigskip


\begin{thebibliography}{99}

\bibitem{A.K.L.S} V. Abramov, R. Kerner, O. Liivapuu, and  S. Shitov, {\it Algebras with ternary law of composition 
and their realization by cubic matrices}, arXiv:0901.2506 

\bibitem{C.RT} R. Campoamor-Stursberg and M. Rausch de Traubenberg, {\it Kinematical superalgebras and 
Lie algebras of order 3}, J.Math.Phys.49:063506 (2008) 

\bibitem{Gne} A. V. Gnedbaye, {\it Op\'{e}rades des alg\`{e}bres $k+1$-aires},
Operads: Proceedings of Renaissance Conferences (Hartford, CT/Luminy, 1995),
Contemp. Math., 202, Amer. Math. Soc., Providence, RI, 1997, pp. 83-113. 

\bibitem{GoM} M. Goze, {\it Alg\`ebres de Lie : classifications, d\'{e}formations et rigidit\'{e},
g\'{e}om\'{e}trie diff\'{e}rentielle.} in {\it Alg\`ebres, dynamique et analyse pour la g\'eom\'etrie}. Edition Ellipse (2009)


\bibitem{G.R1} M. Goze and E. Remm, {\it Lie-admissible algebras and
operads} J. Algebra 273 (2004), no. 1, pp.129-152.

\bibitem{G.R2} M. Goze and E. Remm, {\it Lie admissible coalgebras}, J.
Gen. Lie Theory Appl. 1 (2007), no. 1, pp.19-28.

\bibitem{G.R3} M. Goze and E. Remm, {\it A class of nonassociative algebras},  
Algebra Colloq.  14  (2007),  no. 2, pp. 313-326.

\bibitem{G.R4} E. Remm and M. Goze, {\it On the algebras obtained by
tensor product} arXiv:math/0606105. To appear in Journal of Algebra.

\bibitem{G.R} N. Goze and E. Remm, {\it $n$-ary associative algebras,
cohomology, free algebras and coalgebras}. arXiv:0803.0553 .

\bibitem{H} E. Hoffbeck, {\it A Poincar\'e-Birkhoff-Witt criterion for Koszul
operads},  arXiv:0709.2286.


\bibitem{M.S.S} M. Markl, S. Shnider and J. Stasheff. {\it Operads in algebra,
topology and physics}. Mathematical Surveys and Monographs, 96. American
Mathematical Society, Providence, RI, 2002.

\bibitem{M.RT} G. Moultaka, M. Rausch de Traubenberg and A. Tanasa, {\it A. Cubic supersymmetry 
and abelian gauge invariance},  Internat. J. Modern Phys. A  20  (2005),  no. 25, pp. 5779-5806.

\bibitem{R} E. Remm, {\it On the nonKoszulity of $(2p+1)$-ary partially associative Operads}, 
arXiv:0812.2687.

\bibitem{RT1}  M. Rausch de Traubenberg, {\it Ternary algebras and groups,} 
 Contribution to the 5th International Symposium on Quantum Theory and Symmetries, 
 J.Phys.Conf.Ser.128 (2008).

\bibitem{RT2} M. Rausch de Traubenberg, {\it Cubic extentions of the Poincaré algebra},
Contributed to the XII International Conference on Symmetry Methods in Physics (SYMPHYS-XII), Yerevan, Armenia, July 03-08, 2006
Phys.Atom.Nucl.71 (2008), pp. 1102-1108. 

\end{thebibliography}
\end{document}